\documentclass[a4paper,11pt]{amsart}
\usepackage{amsmath, amsthm, amssymb, amsfonts,amscd}
\usepackage[all]{xy}
\usepackage{graphicx,color}
\usepackage{stmaryrd}

\newtheoremstyle{Style}%
{.5em}{.5em}%
{\it}%
{}%
{\sc}%
{\ {\bf---}}%
{ }%
{}%

\newtheoremstyle{StyleRemarque}%
{.5em}{.5em}%
{\it}%
{}%
{\slshape}%
{.\ }%
{ }%
{}%

\theoremstyle{Style}
\newtheorem{defn}{Definition}[section]
\newtheorem{lem}[defn]{Lemma}
\newtheorem{prop}[defn]{Proposition}
\newtheorem{thm}[defn]{Theorem}

\newtheorem*{thm*}{Theorem}

\theoremstyle{StyleRemarque}

\newtheorem*{rem*}{Remark}

\newtheorem*{ex*}{Example}

\newcommand{\norm}[1]{\left\Vert#1\right\Vert}
\newcommand{\abs}[1]{\left\vert#1\right\vert}

\newcommand{\set}[1]{\left\{#1\right\}}
\newcommand{\Forall}[2]{\forall \, #1 \in #2, \:}

\newcommand{\restric}[1]{\vert_{#1}}
\newcommand{\appli}[3]{#1 \, : \, #2 \To #3}

\newcommand{\rl}{\mathbb R}
\newcommand{\cx}{\mathbb C}
\newcommand{\ir}{\mathbb Z}
\newcommand{\nl}{\mathbb N}
\newcommand{\sph}{\mathbb S}
\newcommand{\bl}{\mathbb B}
\newcommand{\tor}{\mathbb T}
\newcommand{\To}{\longrightarrow}
\newcommand{\OO}{\mathcal O}

\DeclareMathOperator{\tr}{Tr}
\DeclareMathOperator{\ric}{Ric}

\DeclareMathOperator{\vol}{vol}

\DeclareMathOperator{\id}{id}
\DeclareMathOperator{\riem}{Rm}
\DeclareMathOperator{\dvg}{div}

\begin{document}

\title{\bfseries{Rigidity for multi-Taub-NUT metrics.}}
\author{Vincent Minerbe}
\date{\today}

\begin{abstract}
This paper provides a classification result for gravitational instantons with cubic volume growth and cyclic fundamental group at infinity. It proves that a complete hyperk\"ahler manifold asymptotic to a circle fibration over the Euclidean three-space is either the standard $\rl^3 \times \sph^1$ or a	multi-Taub-NUT manifold. In particular, the underlying complex manifold is either $\cx \times \cx/\ir$ or a minimal resolution of a cyclic Kleinian singularity. 
\end{abstract}

\maketitle

\tableofcontents

\section*{Introduction.}

Since the late seventies, there has been considerable interest in the so-called gravitational instantons, namely complete non-compact hyperk\"ahler 
four-manifolds with decaying curvature at infinity. They were introduced by S. Hawking \cite{Haw} as building blocks for his Euclidean quantum 
gravity theory. Beside their natural link with gauge theory, they have also appeared relevant in string theory. Their mathematical beauty,
including the nice twistorial point of view \cite{Bes}, is definitely a good motivation to understand them.

All known examples fall into four families -- ALE, ALF, ALG, ALH -- which differ by their asymptotic geometry. The ALE gravitational instantons are Asymptotically Locally Euclidean, in that their asymptotic geometry is that of $\rl^4$, up to a finite covering. In the other families, the topology outside a compact set and up to finite covering is that of a $\tor^{4-k}$-fibration over $\rl^k$, with $k=3$ for the ALF family, $k=2$ for the ALG family, $k=1$ for the ALH family ; moreover, the geometry is asymptotically adapted to these fibrations. Much more details will be given below, in the ALF case. It is useful to keep in mind that ALE gravitational instantons are exactly those gravitational instantons that have Euclidean volume growth -- $\vol B(x,t) \asymp t^4$ -- \cite{BKN}, while ALF gravitational instantons are characterized by their cubic volume growth -- $\vol B(x,t) \asymp t^3$ \cite{Mi2}.

In 1989, P. Kronheimer classified ALE gravitational instantons \cite{Kr1,Kr2}. The possible topologies are given by the minimal resolutions of the Kleinian singularities $\cx^2/\Gamma$. Here, $\Gamma$ is a finite subgroup of $SU(2)$: cyclic, binary dihedral, tetrahedral, octahedral or icosahedral. For every such manifold $M$, the different hyperk\"ahler structures are parameterized by three classes in $H^{1,1}(M,\cx)$, with a non-degeneracy condition \cite{Kr2}. The simplest situation corresponds to ``$A_{k+1}$'' ALE gravitational instantons, namely $\Gamma=\ir_{k}$, with $k \geq 1$, since they are given by the explicit multi-Eguchi-Hanson metrics. When $k=1$, the manifold is $T^*\cx P^1$ and the metric is also known as Calabi's metric. 

ALF gravitational instantons should be classified in a similar way. What are the examples ? The trivial one is $\cx \times \cx / \ir$, with fundamental group at infinity $\ir$. More sophisticated examples are given by multi-Taub-NUT metrics (cf. section 1), which live exactly on the same manifolds as the multi-Eguchi-Hanson metrics; the fundamental group at infinity is then $\ir_k$. Some more mysterious examples were built recently by S. Cherkis and A. Kapustin, with a dihedral fundamental group at infinity. It is conjectured that these are the only ALF gravitational instantons. 

Indeed, \cite{Mi2} already ensures the topology outside a compact set is like in these examples. Let us recall the precise
statement. We will denote by $\rho$ the distance to some distinguished point. 

\begin{thm}\label{thmvx}\cite{Mi2}
Let $M$ be a four-dimensional complete hyperk\"ahler manifold such that $\riem = \OO(\rho^{-3})$ and 
$\vol B(x,t) \asymp t^3$. Then $M$ is ALF in that there is a compact subset $K$ of $M$ such that $M \backslash K$ is the 
total space of a circle fibration $\pi$ over $\rl^{3}$ or $\rl^{3}/ \{ \pm \id\}$ minus a ball and 
the metric $g$ can be written 
$$
g = \pi^*g_{\rl^3} + \eta^2 + \mathcal{O}(\rho^{-\tau}) \quad \text{for any } \tau <1,
$$  
where $\eta$ is a (local) connection one-form for $\pi$; moreover, the length of the fibers goes to a finite positive limit at infinity. 
\end{thm}

It shows that the topology of the ALF gravitational instanton $M$ can be described outside a compact subset by : $M \backslash K = E \times \rl_+^*$, where $E$ is the total space of a circle fibration over $\sph^2$ or $\rl P^2$. When the base of the circle fibration at infinity is $\sph^2$, we will say that $M$ is \emph{ALF of cyclic type} : in this case, $E$ is either $\sph^2 \times \sph^1$ or $\sph^3/\ir_k$ (where $\ir_k$ is seen as the group of $kth$ roots of $1$, acting by scalar multiplication in $\cx^2$), so the fundamental group of the end is $\ir$ or $\ir_k$, hence the denotation. When the base of the circle fibration at infinity is $\rl P^2$, we say that $M$ is \emph{ALF of dihedral type} : $E$ is then a quotient of $\sph^3$ by a binary dihedral group $D_k$ (of order $4k$), so the fundamental group is dihedral (the trivial bundle is excluded for orientability reasons). In particular, this topological classification for the end rules out any tetrahedral, octahedral or icosahedral fundamental group at infinity. Note that $D_1$ is indeed $\ir_4$, so ``ALF of cyclic type'' means a little bit more than just ``with cyclic fundamental group at infinity''. 

The aim of this paper is to establish a complete classification for ALF gravitation instantons of cyclic type. The precise statement is as follows. 

\begin{thm}
An ALF gravitational instanton of cyclic type is either the flat $\cx \times \cx / \ir$ or a multi-Taub-NUT manifold.
\end{thm}

This classification is up to triholomorphic isometry. In view of Theorem \ref{thmvx}, it means that a complete hyperk\"ahler manifold asymptotic to $\cx \times \cx / \ir$ (resp. a multi-Taub-NUT manifold) is necessarily $\cx \times \cx / \ir$ (resp. a multi-Taub-NUT manifold). In this statement, ``hyperk\"ahler'' cannot be relaxed to ``Ricci-flat''. For instance, there is a non-flat Ricci-flat manifold asymptotic to $\cx \times \cx / \ir$ : the Schwarzschild metric (\cite{Mi3}, for instance). 

The strategy of the proof is as follows. \cite{Mi2} provides an asymptotic model for the manifolds under interest. It involves as key ingredients three functions and a one-form. In the examples, this data extends harmonically to the interior of the manifold and determines the metric. We will prove the existence of such an harmonic extension and then recover the metric from these. In particular, we will extend a ``Killing vector field at infinity'' into a Killing vector field on the whole manifold. 

This work is related to R. Bielawski's paper \cite{Bie}, which (in particular) classifies simply-connected hyperk\"ahler four-manifolds endowed with a trihamiltonian action of $\sph^1$ (see also \cite{KS}). In our context, we are not given a global action of $\sph^1$ but we \emph{build} it, thanks to the Killing vector field mentionned above.   

In a first section, we will describe the multi-Taub-NUT examples, since their properties are crucial for the proof. In a second section, we will start our construction, extending the fibration at infnity into a harmonic map on the whole manifold. In a third section, we will use the complex structures to construct the promised Killing vector field and then finish the proof.

\section{Multi-Taub-NUT metrics.}

Given an integer $k \geq 1$ and $k$ points $a_1,\dots,a_k$ in $\rl^3$, we let $M_{k*}$ be the total space of the principal $\sph^1$ bundle over $\rl^3 \backslash\set{a_1,\dots,a_k}$ whose Chern class integrates to $-1$ over the small two-spheres around each point $a_i$ -- this definition makes sense for these two-spheres form a basis for the homology of $\rl^3 \backslash\set{a_1,\dots,a_k}$ in degree two. Note that the Chern class consequently integrates to $-k$ over the large two-spheres at infinity (Stokes). We need to endow this bundle with a connection. To fix conventions, we recall a connection on a $\sph^1$ bundle is merely a $\sph^1$ invariant one-form $\eta_0$, normalized so that its value on the generator of the action is $1$ (we identify the Lie algebra of $\sph^1$ with $\rl$ and not $i\rl$). It follows that $d\eta_0$ is the pull back of a (``curvature'') two-form $\Omega_0$ on the base, whose cohomology class is the Chern class of the bundle up to a factor $2 \pi$.

Denoting by $x=(x_1,x_2,x_3)$ the coordinates on $\rl^3$, we pick a positive number $m$ and introduce:
$$
V = 1 + \sum_{i=1}^k \frac{2 m}{\abs{x-a_i}}.
$$ 
This a harmonic function so $*_{\rl^3} dV$ is a closed two-form. It follows that $*_{\rl^3} dV$ integrates to $-8m\pi$ over the small two-spheres around each $a_i$, so that $\frac{*_{\rl^3} dV}{4m}$ represents the Chern class of the bundle, hence is the curvature two-form of a connection one-form
$\frac{\eta}{4m}$ (when $k=1$, this is basically the standard contact form on $\sph^3$). The multi-Taub-NUT metrics are then given by the Gibbons-Hawking ansatz:
$$
g = V dx^2 + \frac{1}{V} \, \eta^2 \qquad \text{with } \quad d\eta = *_{\rl^3} dV.
$$
This extends as a complete metric on the manifold $M_k$ obtained by adding one point $p_i$ over each point $a_i$. The $p_i$'s should be thought of as the fixed points of the action of the circle. The ambiguity resulting resulting from the choice of the form $\eta$ only produces isometric metrics (two convenient one-forms $\eta$ only differ by the pull-back of an exact one-form $df$ on the base, which makes them gauge-equivalent: the automorphism $x \mapsto e^{if(x)}\cdot x$ carries one onto the other). 

This metric $g$ turns out be hyperk\"ahler ! The K\"ahler structures are easy to describe. We simply define a complex structure $I$ on $T^*M_{k*}$ by the relations $I dx_1 = \frac{\eta}{V}$ and $I dx_2=dx_3$. One can check that this is indeed K\"ahler \cite{Leb} (hence extends to the whole $M_k$). The other K\"ahler structures are obtained by rotating the roles of $dx_1$, $dx_2$ and $dx_3$. For any of these complex structures, the complex manifold $M_k$ is biholomorphic to the minimal resolution of $\cx^2/\ir_k$ \cite{Leb}. 

At infinity, the curvature decays as $\abs{x}^{-3}$ and the length of the fibers goes to $8m\pi$.

\section{A refined asymptotic model.}

The aim of this section is to improve the asymptotic model provided by \cite{Mi2} for the manifolds we are interested in. In a first step, we recall the rough model from \cite{Mi2}, which is basically a circle fibration at infinity, plus a connection on it. In a second time, we improve the fibration and then the connection. 

\subsection{The rough model.}

Let us give a precise definition for the class of manifolds we are interested in. Given a Riemannian manifold, we denote by $\riem$ the curvature tensor and by $\rho$ the distance function to some distinguished point $o$. We can \emph{define} an ALF gravitational instanton as follows.

\begin{defn}
An ALF gravitational instanton is a complete hyperk\"ahler four-manifold with cubic curvature decay --  $\riem = \OO(\rho^{-3})$ -- and cubic volume growth --
$
\exists \, c>0, \; \forall \, x,\; \forall \, t \geq 1, c^{-1} t^3 \leq \vol B(x,t) \leq c t^3. 
$
\end{defn}
In view of \cite{Mi1}, under the other assumptions, the cubic curvature decay is indeed automatic as soon as the curvature decays faster than quadratically. And it implies the covariant derivatives obey $\nabla^i \riem = \OO(\rho^{-3-i})$ \cite{Mi2}. These facts follow from the Ricci-flatness. 


The paper \cite{Mi2} describes the geometry at infinity of such manifolds. Theorem \ref{thmvx} in the introduction sums up what we need. It ensures the existence of a circle fibration over $\rl^3$ or $\rl^3/\ir_2$ minus a ball at infinity. We \emph{assume} in this paper that the base of this fibration at infinity is $\rl^3$. minus a ball : the ALF gravitational instantons satisfying this property are called ``ALF of cyclic type''. In what follows, we consider an ALF gravitational instanton of cyclic type $(M,g)$.
Let us describe precisely the geometry at infinity, relying on \cite{Mi2}, where every detail is given. 

Basically, $M$ minus a compact subset $K$ is the total space of a circle fibration $\pi$ over $\rl^3$ minus a ball $\bl^3$. Furthermore, this fibration encodes the asymptotic geometry as follows. First, the length of the fibers of $\pi$ goes to some positive and finite value $L_\infty$ at infinity and a $g$-unit $\pi$-vertical vector field $U$ obeys
$$
\nabla^g U = \OO(\rho^{-2}) \quad \text{and} \quad 
\forall\, i \geq 2, \; \nabla^{g,i} U = \OO(\rho^{-i}). 
$$
Second, if we average the metric $g$ into $\widetilde{g}$ along the fibers of $\pi$ (i.e. along the flow of $U$), then
$$
g = \widetilde{g} + \OO(\rho^{-2}) \quad \text{and} \quad 
\Forall{i}{\nl^*} \nabla^{g,i} \widetilde{g} = \OO(\rho^{-1-i}). 
$$
Third, if we push $\widetilde{g}$ down into a metric $\check{g}$ on $\rl^3 \backslash \mathbb{B}^3$, then $\check{g}$ is asymptotically Euclidean of order $\tau$ for any $\tau \in ]0,1[$ in the sense of \cite{BKN}. It implies the existence of $\check{g}$-harmonic coordinates $x_k$ on $\rl^3 \backslash B$ such that 
$$
\check{g} = dx^2 + \OO(\abs{x}^{-\tau}) \quad \text{and} \quad \nabla^{\check{g}}dx_k = \OO(\abs{\check{x}}^{-\tau-1}).
$$

It turns out we can strengthen a little bit the statement about the asymptotic of $\check{g}$. In \cite{Mi2}, we proved a bound on the curvature tensor of $\check{g}$ that ensured this metric was asymptotically Euclidean in $C^{1,\alpha}$ topology, via \cite{BKN}. It turns out we will need a $C^2$ decay of $\check{g}$ to the flat metric. This is only a minor technicality, which we fix now.   

\begin{lem}\label{estimC2}
$\displaystyle{
\nabla^{\check{g},2} dx_k = \OO(\abs{x}^{-\tau-2}).
}$
\end{lem}

\proof
In view of \cite{BKN} (p. 314-315), if we prove $\nabla^{\check{g}} \riem_{\check{g}} = \OO(\abs{x}^{-4})$, then the asymptotically Euclidean behaviour is
true with order $\tau$ in $C^{2,\alpha}$ topology and we are done. To prove this estimate, we 
choose exponential coordinates at the running point on the base and lift the coordinate vector fields $\partial_i$ into vector fields $X_i$ that are $\widetilde{g}$-orthogonal to the fibers of $\pi$. O'Neill's formula (\cite{Bes}) expresses the quantity $\check{g}(\riem_{\check{g}} (\partial_i,\partial_j) \partial_k, \partial_l )$ as $\widetilde{g}(\riem_{\widetilde{g}} (X_i,X_j) X_k, X_l )$ plus a linear combination of terms like $\widetilde{g}([X_i,X_j],U) \widetilde{g}([X_k,X_l],U)$. The formula can be differentiated to get 
$$
\abs{\nabla^{\check{g}}\riem_{\check{g}}}
\leq c \abs{\nabla^{\widetilde{g}}\riem_{\widetilde{g}}}
+ c \abs{\nabla^{\widetilde{g}}U} \left(\abs{\riem_{\widetilde{g}}}
+ \abs{\nabla^{\widetilde{g},2}U} + \abs{\nabla^{\widetilde{g}}U}^2 
\right).
$$
Using the estimates recalled above, we find $\nabla^{\widetilde{g},i}\riem_{\widetilde{g}} = \OO(\rho^{-3-i})$,
$\nabla^{\widetilde{g}}U = \OO(\rho^{-2})$ and $\nabla^{\widetilde{g},2}U = \OO(\rho^{-2})$,
hence the result.
\endproof

\subsection{A best fibration at infinity.}

The estimates above make it possible to find $g$-harmonic functions that approach the (pullback of the) functions $x_k$ at infinity in $C^1$ topology, with appropriate $C^\infty$ estimates ; a ``best'' fibration $\pi$ will stem from these harmonic functions. Beware 
we will often use the same notation for a function on the base $\rl^3$ and its pullback by the fibration. 

\begin{lem}\label{estimnewfib}
For every index $k$, one can find a $g$-harmonic function $\underline{x}_k$ on $M$ such that
for any $\epsilon >0$,
$
\underline{x}_k = x_k + \OO(\rho^{\epsilon}) \quad \text{and} \quad d\underline{x}_k = dx_k + \OO(\rho^{\epsilon-1}),
$
with moreover:
$\displaystyle{
\forall i\geq 2, \; \nabla^{g,i} \underline{x}_k = \OO(\rho^{\epsilon-i}).
}$
\end{lem}

\proof
Let us extend $x_k$ as a smooth function on the whole $M$ and observe 
$$
\abs{\Delta_g x_k - \Delta_{\widetilde{g}} x_k} 
\leq c \abs{g-\widetilde{g}} \abs{\nabla^{\widetilde{g}} dx_k} + c \abs{\nabla^{g} - \nabla^{\widetilde{g}}} \abs{dx_k} \leq c \rho^{-2}. 
$$
Since the functions $x_k$ are harmonic with respect to $\check{g}$ or $\widetilde{g}$, we deduce $\Delta_g x_k = \OO(\rho^{-2})$.  This lies in $\rho^{\delta-2} L^2$ for any $\delta > \frac{3}{2}$, so we can apply the analysis of \cite{Mi3} to find a solution $u_k$ for the equation 
$\Delta_g u_k = - \Delta_g x_k$, with $\nabla^i u_k \in \rho^{\delta-i} L^2$, $0\leq i \leq 2$. As explained in the appendix of \cite{Mi2}, a Moser iteration yields
$$
\norm{u_k}_{L^\infty(A_R)} \leq c R^{-\frac{3}{2}} \norm{u_k}_{L^2(A'_R)} 
+ c R^2 \norm{\Delta_g u_k}_{L^\infty(A'_R)},
$$
where $A_R = \set{R\leq \rho \leq 2R}$ and $A'_R = \set{R/2\leq \rho \leq 4R}$. Since $u_k$ is in $\rho^{\frac{3}{2}+\epsilon} L^2$ and $\Delta_g u_k = - \Delta_g x_k = \OO(\rho^{-2})$, we get $u_k = \OO(\rho^\epsilon)$ for any positive $\epsilon$. 
Since $\ric_g = 0$, the Hodge Laplacian and the Bochner Laplacian (which we denote by $\Delta_g$ on the whole tensor algebra) coincide on one-forms. With lemma \ref{estimC2}, we can apply the same argument to $du_k$, with this Laplacian (cf. \cite{Mi2}) and find $du_k = \OO(\rho^{\epsilon-1})$. As a result, the function $\underline{x}_k := x_k + u_k$ is $g$-harmonic, with    
$$
\underline{x}_k = x_k + \OO(\rho^{\epsilon}) \quad \text{and} \quad d\underline{x}_k = dx_k + \OO(\rho^{\epsilon-1}).
$$
The equation $\Delta_g d\underline{x}_k = 0$ can be used together with
$
[\Delta_g, \nabla^g] = \riem_g \varoast \nabla^g + \nabla^g \riem_g \varoast
$
(here, $\varoast$ denotes any bilinear pairing depending only on $g$) to obtain:
$$
\Delta_g \nabla^{g,i} d\underline{x}_k 
= \sum_{j=0}^{i} \nabla^{g,j} \riem_g \varoast \nabla^{g,i-j} d\underline{x}_k.
$$ 
Since $\nabla^{g,j} \riem_g = \OO(\rho^{-3-j})$, the estimates follow from an induction argument based on the following inequality (Moser iteration and cutoff argument, cf. \cite{Mi2}):
\begin{eqnarray*}
& & \norm{\nabla^{g,i} d\underline{x}_k}_{L^\infty(A_R)} + R^{-\frac{1}{2}} \norm{\nabla^{g,i+1} d\underline{x}_k}_{L^2(A_R)} \\
&\leq& c R^{-\frac{3}{2}} \norm{\nabla^{g,i} d\underline{x}_k}_{L^2(A'_R)} 
+ c R^{\frac{1}{2}} \norm{\Delta_g \nabla^{g,i} d\underline{x}_k}_{L^2(A'_R)}.
\end{eqnarray*}
\endproof

Enlarging $K$ if necessary, these $g$-harmonic functions $\underline{x}_1$, $\underline{x}_2$, $\underline{x}_3$ provide a new $\sph^1$-fibration 
$
\appli
{ \underline{\pi} = (\underline{x}_1,\underline{x}_2,\underline{x}_3) }
{M \backslash K }
{\rl^3 \backslash B}
$. As a consequence of (\ref{estimnewfib}), vector fields $X$ that are $g$-orthogonal to the fibers of $\underline{\pi}$ obey
\begin{equation}\label{quasiiso}
\frac{\abs{d\underline{\pi}(X)}_{g_{\rl^3}}}{\abs{X}_g} = 1 + \OO(\rho^{\epsilon-1})
\end{equation}
and $\underline{\pi}$ satisfies the following estimates  
\begin{equation}\label{estimsup}
\forall i\geq 2, \; \nabla^{g,i} \underline{\pi} = \OO(\rho^{\epsilon-i}),
\end{equation}
with respect to the metric $g$ on $M\backslash K$ and the Euclidean metric on the base. 
Let us denote by $\underline{U}$ a $g$-unit $\underline{\pi}$-vertical vector field. Differentiating
the relation $d\underline{\pi} (\underline{U}) =0$ and using (\ref{quasiiso}), (\ref{estimsup}), we get as in \cite{Mi2}:
\begin{equation}\label{estimV}
\Forall{i}{\nl^*}  \nabla^{g,i} \underline{U} = \OO(\rho^{\epsilon-1-i}).
\end{equation}
We also introduce the function $\underline{L}$ assigning to each point $p$ the length of the $\underline{\pi}$-fiber
through $p$ and the flow $\underline{\psi}^t$ of $\underline{U}$. With $\underline{\psi}^L = \id$, one finds
$
d\underline{L} = \left( g - \underline{\psi}^{\underline{L}*} g \right)(\underline{U},.).
$  
Since the Lie derivative of $g$ along $\underline{U}$ is twice the symmetrization of $\nabla^g \underline{U}$, we have
$$
\abs{g - \underline{\psi}^{t*}g} \leq c t \rho^{\epsilon-2}
\quad \text{and therefore} \quad
\abs{d\underline{L}} \leq c \underline{L} \rho^{\epsilon-2}.
$$
So $d\log \underline{L}= \OO(\rho^{\epsilon-2})$. From Cauchy's criterion, we deduce $\underline{L}$ goes to 
a limit $\underline{L}_\infty$ at infinity (and $\underline{L}_\infty$ = $L_\infty$, indeed).
Using the arguments of \cite{Mi2}, we can then estimate the metric $\widetilde{\underline{g}}$ obtained by averaging 
$g$ along the flow of $\underline{U}$ by
\begin{equation}\label{symclose}
\Forall{i}{\nl} \nabla^{g,i}( \widetilde{\underline{g}} - g) = \OO(\rho^{\epsilon-2-i})
\end{equation}
and the derivatives of $\underline{L}$ by
$$
\Forall{i}{\nl^*} \nabla^{g,i}(\underline{L}) = \OO(\rho^{\epsilon-1-i}).
$$
We introduce the vector field $\underline{T} := \frac{\underline{L}}{\underline{L}_\infty} \underline{U}.
$ on $M\backslash K$ and see $\frac{\underline{L}_\infty}{2\pi} \underline{T}$ as the infinitesimal generator of a $\sph^1$ action, which 
makes $\underline{\pi}$ into a principal $\sph^1$-fibration. The one-form $
\underline{\eta} := \frac{\widetilde{\underline{g}}(\underline{T},.)}{\widetilde{\underline{g}}(\underline{T},\underline{T})} 
$
is $\sph^1$ invariant and satisfies $\underline{\eta}(\underline{T})=1$. It is $\frac{\underline{L}_\infty}{2\pi}$ times a connection on the $\sph^1$-bundle. The metric $\widetilde{\underline{g}}$ induces a metric $\check{\underline{g}}$ on the base and we have
$$
\Forall{i}{\nl}  \nabla^{\check{\underline{g}},i} \left( \check{\underline{g}} - d\underline{x}^2 \right)
= \OO(\abs{\check{\underline{x}}^{\epsilon-1-i}}). 
$$
The outcome of all this is a ``model'' metric  
$
\underline{h} := d\underline{x}^2 + \underline{\eta}^2 
$
such that $d\underline{\eta} = \underline{\pi}^*\underline{\Omega}$ with
$$
\Forall{i}{\nl} \nabla^{\underline{h},i} (g-\underline{h}) = \OO(\rho^{\epsilon-1-i})
\quad \text{and} \quad 
D^{i}(\underline{\Omega}) = \OO(\abs{\underline{x}}^{\epsilon-2-i}).
$$
For the sake of simplicity, we will forget the underlining and denote $\abs{x}$ by $r$. The letter $\epsilon$ will refer to any small positive number and we basically use the convention $"\epsilon=2\epsilon"$ (finitely many times, hopefully).

What have we gained ? We chose to sacrifice an $\epsilon$ in the lower order estimates of $\pi$ and $U$ for two benefits. The first one is technical : we are provided with better estimates for the higher order derivatives. The second one is essential: the fibration at infinity extends as a harmonic function on the whole manifold, as in the examples. It is encouraging to notice that such a function $\pi$ is unique up to the addition of a constant: this is certainly the good object to look at!  

\section{Using the hyperk\"ahler structure}

\subsection{A Killing vector field.}

We can rely on the hyperk\"ahler structure $(g,I,J,K)$ to build a Killing vertical vector field. To fix ideas, let us proceed to some normalization. Since $d (dx_k(IT))=\OO(r^{\epsilon-2})$, the functions $dx_k(IT)$ go to some constants at infinity (Cauchy criterion). Indeed, since moreover  $\eta(IT) =\OO(r^{\epsilon-1})$, we can rotate the coordinates $x_k$ so that $dx_k(IT)=\delta_{1k}+ \OO(r^{\epsilon-1})$. Similarly, up to a second rotation (in the plane $x_1=0$ in $\rl^3$), we can assume $dx_k(JT)=\delta_{2k} + \OO(r^{\epsilon-1})$, and consequently $dx_k(KT)=\delta_{3k} + \OO(r^{\epsilon-1})$. 

The following proposition is a key step. Its proof again uses the fact that a hyperk\"ahler metric has vanishing Ricci curvature, so that the Hodge Laplacian on $1$-forms coincides with the Bochner Laplacian. 

\begin{prop}
There is a unique $g$-Killing vector field $W$ such that
$$
\iota_W \omega_I = -dx_1, \quad \iota_W \omega_J = -dx_2, \quad \iota_W \omega_K = -dx_3.
$$
It is $\pi$-vertical, preserves the complex structures $I$, $J$, $K$ and obeys the estimate
$$
\Forall{k}{\nl} \nabla^{g,k} (W-T) = \OO(r^{\epsilon-1-k}).
$$
\end{prop}

\proof
One can define three vector fields $W_1$, $W_2$, $W_3$ by the relations
$\iota_{W_1} \omega_I = -dx_1$, $\iota_{W_2} \omega_J = -dx_2$ and $\iota_{W_3} \omega_K = -dx_3$.
Since the one-forms $dx_l$ are $g$-harmonic and the K\"ahler forms are parallel, the vector fields 
$W_l$ are $g$-harmonic. Now, for $l=1,2,3$, our previous estimates ensure 
$$
\nabla^{k} (W_l-T) = \OO(r^{\epsilon-1-k}),
$$
for every \emph{positive} integer $k$. Our choice of coordinates in $\rl^3$ moreover implies that $W_l -T$ goes to zero at infinity. One can therefore integrate the estimate for $k=1$ into the corresponding estimate for $k=0$. In particular, the difference $X$ between any two of these vector fields $W_l$ is harmonic and goes to zero at infinity. It follows that the function $\abs{X}^2$ goes to zero at infinity and satisfies :
$$
\Delta \abs{X}^2 = 2 g(\Delta X,X) - 2 \abs{\nabla X}^2 \leq 0.
$$
From the maximum principle, we deduce $\abs{X}^2=0$, namely $X=0$. We conclude: $W_1=W_2=W_3=:W$. By definition, we 
have $dx_1(W) = - \iota_W \omega_I(W)=0$ as well as $dx_2(W)=dx_3(W)=0$ for the same reason. So $W$ is vertical. Since $\omega_I$ and $\iota_W \omega_I=-dx_1$ are closed, Cartan's magic formula yields $L_W \omega_I = 0$; similarly, $L_W \omega_J = L_W \omega_K = 0$. Finally, since $\omega_I$ is parallel, we have: $\iota_{\nabla W} \omega_I = - \nabla dx_1$. The right-hand side is a symmetric bilinear form, so for any two vector fields $X$, $Y$, we are given the identity: $g(\nabla^g_X W,Y) = -g(\nabla^g_{I Y} W, I X)$. We can of course play the same game with $K$ and then $J$, so as to find
$$
g(\nabla^g_X W,Y) = g(\nabla^g_{K I X} W, K I Y) = -g(\nabla^g_{J K I Y} W, J K I X) = - g(\nabla^g_{Y} W, X).
$$
The tensor $\nabla^g W$ is therefore skew-symmetric, which exactly means $W$ is Killing. Since the K\"ahler forms and the metric are preserved by $W$, the complex structures are also preserved.
\endproof

Let us set $\alpha:= g(W,.)$. The definition of $W$ means $dx_1= I\alpha$, $dx_2= J\alpha$
and $dx_3= K\alpha$. In particular, the covectors $dx_k$ are everywhere $g$-orthogonal and have the same $g$-norm as $\alpha$ -- notice these covectors thus vanish only all at the same time. With $V:= \abs{W}^{-2}$, we might keep in mind the following formulas, on the open set $M_*$ where $W$ does not vanish: 
$$
g = V (dx^2 + \alpha^2) \quad \text{and} \quad \omega_I = V \left( dx_1 \wedge \alpha + dx_2 \wedge dx_3 \right).
$$

\subsection{The map $\pi$.}

\begin{lem}
The set $M \backslash M_*$ is finite : $M \backslash M_* = \set{p_1,\dots,p_k}$.  
\end{lem}

\proof
$M \backslash M_*$ is the place where the vector field $W$ vanishes. Pick a point $p$ such that $W_p=0$. The flow $\phi^t$ of $W$ 
preserves the hyperk\"ahler structure so the differential $T_p \phi^t$ acts on $T_p M$ by a transformation in $SU(2)$, reading 
$
\left(
\begin{array}{cc}
e^{i\lambda t} & 0 \\
0 & e^{-i\lambda t}
\end{array}
\right)
$
in some basis, for some $\lambda \in \rl^*$ ; this is indeed a special case of \cite{GH}. It follows that $p$ is an isolated fixed point
of the flow. This ensures $M \backslash M_*$ is discrete. Since $V$ goes to $1$ at infinity, $M \backslash M_*$ is compact. It is 
therefore a finite set.   
\endproof

\begin{lem}
$\appli{\pi}{M}{\rl^3}$ is onto. 
\end{lem}

\proof
Since $\pi\restric{M_*}$ is a submersion, the set $\pi(M_*)$ is open in $\rl^3$. Let 
$y$ be a point of its boundary and let $y_n = \pi (x_n)$ denote a sequence in $\pi(M_*)$ 
and converging to $y$. Since $\pi$ is asymptotic to a circle fibration with fibers of bounded
length, $\pi$ is proper, so we may assume $x_n$ goes to some point $x$ in $M$, with $\pi(x)=y$.
Since $\pi(M_*)$ is open, $x$ cannot belong to $M_*$ : $x$ is one of the $p_i$'s. To sum up, $\pi(M_*)$ is an open set of $\rl^3$ whose boundary is contained in $\set{\pi(p_1),\dots,\pi(p_k)}$. In particular, $\pi(M)$ is dense in $\rl^3$. It is also closed, for $\pi$ is proper : $\pi(M) = \rl^3$.  
\endproof

Let us set $a_i = \pi(p_i)$.

\begin{lem}
The map $\appli{\pi}{M_*}{\rl^3 \backslash \set{a_1,\dots,a_k}}$ is a circle fibration. 
\end{lem}

\proof
This submersion $\appli{\pi}{M_*}{\rl^3 \backslash \set{a_1,\dots,a_k}}$ is surjective and proper (it can be proved exactly as in the lemma above), hence a fibration by Ehresmann's theorem. Since it is a circle fibration outside a compact set, it is a circle fibration.
\endproof

\begin{rem*}
$\pi$ can be interpreted as the hyperk\"ahler moment map for the Hamiltonian action given by the flow of $W$. This kind of situation is studied in \cite{Bie}, which inspired us.  
\end{rem*}

\subsection{The connection.}

Let us \emph{redefine} the one-form $\eta := V \alpha$, so that $\eta(W)=1$, $L_W\eta=0$ and thus $\iota_W d\eta=0$ ; it follows that $d\eta=\pi^*\Omega$ for some two-form $\Omega$ on the base. We wish to interpret $\eta$ as a connection form (up to some normalization).  

The fibers of $\pi$ are also the orbits of $W$, so that the flow of $W$ is periodic. In $M_*$, the period $P_x$ of the orbit $\pi^{-1}(x)$ is equal to the integral $\int_{\pi^{-1}(x)} \eta$. Given nearby points $x$ and $y$, it follows from Stokes theorem that the difference $P_x -P_y$ is the integral of $d\eta$ over the cylinder $\pi^{-1}([x,y])$, which vanishes because $d\eta$ is basic! The period is therefore constant and can be computed at infinity : it is equal to the length $L_\infty$ of the fibers at infinity. In particular, the vector field $\frac{L_\infty}{2\pi} W$ is the infinitesimal generator of an action of $\sph^1$, making $\pi\restric{M\backslash K}$ a principal $\sph^1$ bundle over $\rl^3$ minus a ball. Furthermore, $\frac{2\pi}{L_\infty} \eta$ is a connection one-form on this $\sph^1$-bundle and its curvature can be computed through the following lemma.

\begin{lem}\label{courbure} 
$\displaystyle{\Omega = *_{\rl^3} dV}$.
\end{lem}

\proof
The K\"ahler form $\omega_I =  dx_1 \wedge \alpha + V dx_2 \wedge dx_3$ is closed, hence the relation
$
dx_1\wedge d\eta = dV \wedge dx_2 \wedge dx_3,
$
which means: $\Omega(\partial_2,\partial_3)= \partial_1 V$. The other K\"ahler forms provide the remaining components.
\endproof

So we know $\Omega$ as soon as we know $V$.

\begin{lem}\label{teteF}
There are positive numbers $m_i$ such that 
$\displaystyle{
V = 1 + \sum_{i=1}^k \frac{2m_i}{\abs{x - a_i}}.
}$
\end{lem}

\proof
Lemma \ref{courbure} implies $V$ is $g_{\rl^3}$-harmonic outside the $a_i$'s. Moreover, it is positive. A classical 
result \cite{ABR} then ensures that around each singularity $a_i$, it is the sum of the function $\frac{2m_i}{\abs{x - a_i}}$, $m_i >0$, and of a smooth harmonic function. Globally, this means  
$
V = \varphi + \sum_{i=1}^k \frac{2m_i}{\abs{x - a_i}}
$
for some smooth harmonic function $\varphi$ on $\rl^3$. The asymptotic of the metric implies $\varphi$ goes to $1$ at infinity, so that $\varphi = 1$.  \endproof

We can then identify $\Omega$. If $d\omega_i$ is the volume form of the unit sphere around $a_i$, we have : 
$$
\Omega = -\sum_{i=1}^k 2 m_i d\omega_i.
$$
As recalled in section 1, this data determines the connection up to gauge (which is enough for a classification up to isometry). 

Observe the topology of the circle bundle determines the cohomology class of $\Omega$, seen as $\frac{L_\infty}{2\pi}$ times the curvature of a connection. For large $R$, this means $\frac{1}{L_\infty} \int_{r=R} \Omega = c_1^\infty$, where $c_1^\infty$ is the Chern number of the fibration over the large spheres. A relation follows: 
$$
-8\pi  \sum_{i=1}^k m_i = L_\infty c_1^\infty.
$$
Let us now look at the circle bundle induced on a small sphere near $a_i$. It has Chern number $c_1^i =\pm 1$ for there is no orbifold singularity on $M$ and it can be computed by an integral of $\Omega$ as above:
$$
-8\pi  m_i = L_\infty c_1^i.
$$
This has two consequences. First, $c_1^i= -1$ because $m_i$ is nonnegative. Second, $8\pi m_i = L_\infty$ ! The parameters $m_i$ are necessarily all equal, depending only on the length of the fibers at infinity (a similar remark can be found \cite{GH}).

Observe also that the number $k$ of singularities is given by the Chern number at infinity : $k=-c_1^\infty$ ; the fundamental group at infinity is simply $\ir_k$. The topology is completely determined by the parameter $k$. 

\begin{rem*}
It is interesting to relate the asymptotic of $V$ to the mass $m$ defined in this paper. It is a nonnegative Riemannian invariant, given by
$$
m := -\frac{1}{12\pi L_\infty} \lim_{R \To \infty} \int_{\partial B_R} *_h  \; \Big{(} \dvg_h g + d \tr_h g -\frac{1}{2} d \; g(W,W) \Big{)},
$$
where $h=dx^2+\eta^2$ in our context (this definition indeed differs by a factor $3$ from that in \cite{Mi3}). Here, a slight computation (cf. \cite{Mi3}) ensures this mass is $m = \sum_{i=1}^k m_i$, hence $-8\pi  m = L_\infty c_1$. The classification result when $k=0$ is an immediate application of \cite{Mi3} for the mass $m$ vanishes iff $(M,g)$ is isometric to the standard $\rl^3 \times \sph^1$.  
\end{rem*}

\subsection{The classification.}

We can conclude by the

\begin{thm}
When $k=0$, $(M,g)$ is isometric to the standard $\rl^3 \times \sph^1$, with no holonomy and circles of length $L_\infty$. When $k\geq 1$, $(M,g)$ is isometric to $M_k$ endowed with the multi-Taub-NUT metric given by
$$
g= V dx^2 + \frac{1}{V} \eta^2,
$$
where $V = 1 + \sum_{i=1}^k \frac{2m}{x-a_i}$ and $d\eta = *_{\rl^3} dV$. The positive parameter $m$ is the mass of $(M,g)$ and it is $\frac{k}{8\pi}$ times the length of the fibers at infinity. 
\end{thm}

The holomorphic classification follows from the explicit formulas for the K\"ahler forms. 



\end{document}